\documentclass[aop,MSNbibl,seceqn,dvips]{arximspdf}
\usepackage{mathrsfs}

\doi{10.1214/12-AOP815} 
\volume{41}
\issue{4}
\pubyear{2013}
\firstpage{2709}
\lastpage{2723}

\makeatletter
\newtheorem{prop}{Proposition}[section]

\newtheorem{lemma}[prop]{Lemma}
\newproclaim{rem}[prop]{Remark}
\newproclaim{rems}[prop]{Remarks}
\newtheorem{thmm}[prop]{Theorem}
\newproclaim{defi}[prop]{Definition}
\newproclaim{example}[prop]{Example}
\newproclaim{examples}[prop]{Examples}

\newcommand{\R}{\mathbb{R}}

\newcommand{\cont}[1]{\displaystyle{\mathop{\frown}^{#1}}}

\makeatother

\begin{document}
\begin{frontmatter}

\title{Poisson approximations on the free Wigner chaos}
\runtitle{Poisson approximations on the free Wigner chaos}

\begin{aug}
\author[A]{\fnms{Ivan} \snm{Nourdin}\ead[label=e1]{inourdin@gmail.com}}
\and
\author[B]{\fnms{Giovanni} \snm{Peccati}\corref{}\ead[label=e2]{giovanni.peccati@gmail.com}}
\runauthor{I. Nourdin and G. Peccati}
\affiliation{Universit\' e De Lorraine and Universit\'{e} du Luxembourg}
\address[A]{Institut \'Elie Cartan\\
Universit\' e De Lorraine\\
BP 70239, 54506 Vandoeuvre-l\`es-Nancy\\
France\\
\printead{e1}} 
\address[B]{Facult\'{e} des Sciences\\
de la Technologie et de la Communication\\
UR en Math\'{e}matiques. 6\\
rue Richard Coudenhove-Kalergi\\
L-1359 Luxembourg\\
Luxembourg\\
\printead{e2}}
\end{aug}

\received{\smonth{5} \syear{2011}}
\revised{\smonth{10} \syear{2012}}

%
\begin{abstract}
We prove that an adequately rescaled sequence $\{F_n\}$ of self-adjoint
operators, living inside a fixed
free Wigner chaos of even order, converges in distribution to a
centered free Poisson random variable with rate $\lambda>0$ if and
only if
$\varphi(F_n^4)- 2\varphi(F_n^3)\rightarrow2\lambda^2-\lambda$ (where
$\varphi$ is the relevant tracial state). This extends to a free
setting some recent limit theorems by Nourdin and Peccati [\textit{Ann. Probab.} \textbf{37} (2009) 1412--1426]
and provides a noncentral
counterpart to a result by Kemp et al. [\textit{Ann. Probab.} \textbf{40} (2012) 1577--1635].
As a by-product of our findings, we show that
Wigner chaoses of order strictly greater than 2 do not contain nonzero
free Poisson random variables. Our techniques involve the so-called
``Riordan numbers,'' counting noncrossing partitions without singletons.
\end{abstract}

%

\end{frontmatter}
%
\section{Introduction}\label{sec1}

\subsection{Overview}\label{sec1.1}

Let $W$ be a standard Brownian motion on $\R_+$, and let $q\geq1$ be
an integer. For every deterministic symmetric function
$f \in L^2(\R_+^q)$, we denote by $I_q^W(f)$ the multiple stochastic
Wiener--It\^o integral of $f$ with respect to $W$. Random variables of
the form $I_q^W(f)$ compose the so-called $q$th \textit{Wiener chaos}
associated with $W$. The concept of Wiener chaos roughly represents an
infinite-dimensional analogous to Hermite polynomials for the
one-dimensional Gaussian distribution; see, for example,~\cite{PecTaq}
for an introduction to this topic.

The following two results, proved, respectively, in~\cite{NP} and \cite
{NPchi}, provide an exhaustive characterization of normal and Gamma
approximations on Wiener chaos. As in~\cite{NPchi}, we denote by
$F(\nu
)$ a centered random variable with the law of $2G(\nu/2)-\nu$, where
$G(\nu/2)$ has a Gamma distribution with parameter $\nu/2$ [if $\nu
\geq
1$ is an integer, then $F(\nu)$ has a centered $\chi^2$ distribution
with $\nu$ degrees of freedom].

\begin{thmm}\label{TNPEC} \textup{(A)} 
Let $N\sim\mathscr{N}(0,1)$, fix $q\geq2$ and let $I^W_q(f_n)$ be a
sequence of multiple stochastic integrals with respect to the standard
Brownian motion $W$, where each $f_n$ is a symmetric element of $L^2(\R_+^q)$ such that $E[I_q^W(f_n)^2] = q!\|f_n\|^2_{L^2(\R_+^{q})} = 1$.
Then, the following two assertions are equivalent, as $n\to\infty$:
\begin{longlist}[(ii)]
\item[(i)] $I^W_q(f_n)$ converges in distribution to $N$;
\item[(ii)] $E[I^W_q(f_n)^4]\to E[N^4] = 3$.
\end{longlist}

\textup{(B)} 
Fix $\nu>0$, and let $F(\nu)$ have the centered Gamma distribution
described above. Let $q\geq2$ be an even integer, and
let $I^W_q(f_n)$ be a sequence of multiple stochastic integrals, where
each $f_n$ is symmetric and verifies $E[I_q^W(f_n)^2] =E[F(\nu
)^2]=2\nu$.
Then the following two assertions are equivalent, as $n\to\infty$:
\begin{longlist}[(ii)]
\item[(i)] $I^W_q(f_n)$ converges in distribution to $F(\nu)$;
\item[(ii)] $E[I_q(f_n)^4] -12E[I_q(f_n)^3] \to E[F(\nu)^4]
-12E[F(\nu)^3] = 12\nu^2 - 48\nu$.
\end{longlist}
\end{thmm}

The results stated in Theorem~\ref{TNPEC} provide a drastic
simplification of the so-called \textit{method of moments} for probabilistic
approximations, and have triggered a huge amount of applications and
generalizations, involving, for example, Stein's method, Malliavin calculus,
power variations of Gaussian processes, Edgeworth expansions, random
matrices and universality results. See~\cite{NPptrf,NPsurvey}, as
well as
the monograph~\cite{NPbook}, for an overview of the most important
developments. See~\cite{WWW} for a constantly updated web
resource, with links to all available papers on the subject.

In~\cite{KNPS}, together with Kemp and Speicher, we proved an analogue
of part (A) of Theorem~\ref{TNPEC} in the framework of free probability
(and free Brownian motion). Let $(\mathscr{A}, \varphi)$ be a free
probability space, and let $\{S(t) \dvtx t\geq0\}$ be a free Brownian motion
defined therein; see Section~\ref{swigner} for details. As shown
in~\cite{bianespeicher}, one can define multiple integrals of the type
$I_q^S(f)$, where $f$ is a square-integrable complex kernel [to
simplify the notation, throughout the paper we shall drop the suffixes
$q,S$, and write simply $I(f) = I_q^S(f)$]. Random variables of the
type $I(f)$ compose the so-called \textit{Wigner chaos} associated with
$S$, playing in free stochastic analysis a
role analogous to that of the classical Gaussian Wiener chaos; see, for
instance,~\cite{bianespeicher}, where Wigner chaoses are used to
develop a free version of the Malliavin calculus of variations. The
following statement is the main result of~\cite{KNPS}.

\begin{thmm}\label{TKNPS} 
Let $s$ be a centered semicircular random variable with unit variance
[see Definition~\ref{DZ}\textup{(i)}], fix an integer $q\geq2$ and
let $I(f_n)$ be a sequence of multiple integrals of order $q$ with
respect to the free Brownian\vadjust{\goodbreak} motion $S$, where each $f_n$ is a mirror symmetric
(see Section~\ref{swigner}) element of $L^2(\R_+^q)$ such that
$\varphi
[I_q(f_n)^2] = \|f_n\|^2_{L^2(\R_+^{q})} = 1$. Then the following two
assertions are equivalent, as $n\to\infty$:
\begin{longlist}[(ii)]
\item[(i)] $I(f_n)$ converges in distribution to $s$;
\item[(ii)] $\varphi[I(f_n)^4]\to\varphi[s^4] = 2$.
\end{longlist}
\end{thmm}

The principal aim of this paper is to prove a free analogy of part (B) of
Theorem~\ref{TNPEC}. As explained in Section~\ref{SPoisson}, and
somewhat counterintuitively, the free analogy of Gamma random variables
is given by free Poisson random variables; see Definition~\ref{DZ}(ii).\vspace*{-6pt}

\begin{rem}
(i) The counterintuitive nature of the correspondence between
the free Gamma and the free Poisson distribution appears most
prominently when one considers a free Poisson random variable $Z(p)$
with integer parameter $p\in\{1,2,\ldots\}$. In this case, the law of
$Z(p)$ is both equal to the law of the sum of $p$ freely independent
squared semicircular random variables (a proof of this fact is provided
in Proposition~\ref{Po}), and to the limit of some appropriate free
convolution of Bernoulli distributions; see~\cite{nicaspeicher}, Proposition~12.11.
This correspondence has of course no analogous in
classical probability. As explained in~\cite{nicaspeicher}, Remark~12.14, such a phenomenon is one of the many
manifestations of the specific algebraic structure of the lattice
$NC(n)$ of all noncrossing partitions of the set $[n] =\{1,\ldots,n\}$
($n=1,2,\ldots$), in terms of which free cumulants are expressed. In
particular, the lattice $NC(n)$ is \textit{self-dual}, with the duality
implemented by the so-called \textit{Kreweras complementation}; see
\cite{nicaspeicher}, page 147. No such self-dual structure exists for the
lattice $P(n)$ of all partition of $[n]$, playing a role analogous to
$NC(n)$ in the computation of classical cumulants (see, e.g.,
\cite{PecTaq},
Chapter 3), and it is exactly this lack of additional
symmetry that explains the combinatorial difference between Gamma and
Poisson distributions in a classical framework.\vspace*{-6pt}
\begin{longlist}[(ii)]
\item[(ii)] The free Poisson law is also known as the
\textit{Marchenko--Pastur distribution}, arising in random matrix
theory as the limit of the eigenvalue distribution of large sample
covariance matrices; see, for example, Bai and Silverstein~\cite{bs},
Chapter~3, Hiai and Petz~\cite{HD}, pages 101--103 and~130, and the
references therein.
\end{longlist}
\end{rem}

The following statement is the main achievement of the present work.
%
\begin{thmm}\label{mainthm} Let $q\geq2$ be an even integer. Let
$Z(\lambda)$ have a centered free Poisson distribution with rate
$\lambda>0$. Let $I(f_n)$ be a sequence of multiple integrals of order
$q$ with respect to the free Brownian motion $S$, where each $f_n$ is
a mirror symmetric element of $L^2(\R_+^q)$ such that $\varphi
[I_q(f_n)^2] = \|f_n\|^2_{L^2(\R_+^{q})}= \varphi[Z(\lambda)^2] =
\lambda$.\vadjust{\goodbreak}
Then the following two assertions are equivalent, as $n\to\infty$:
\begin{longlist}[(ii)]
\item[(i)] $I(f_n)$ converges in distribution to $Z(\lambda)$;
\item[(ii)] $\varphi[I(f_n)^4] - 2\varphi[I(f_n)^3] \to
\varphi
[Z(\lambda)^4]- 2\varphi[Z(\lambda)^3] = 2\lambda^2-\lambda$.
\end{longlist}
\end{thmm}

One should note that the techniques involved in our proofs are
different from those adopted in the previously quoted references,
as they are based on a direct enumeration of contractions. These
contractions emerge when iteratively applying product formulas for multiple
Wigner integrals; see also~\cite{NouYet}. One crucial point is that the
moments of a free Poisson random variable can be expressed in terms of
the so-called \textit{Riordan numbers}, counting the number of
noncrossing partitions without singletons; see, for example,~\cite{bernhart}.
We also stress that one cannot expect to have convergence to a nonzero
free Poisson inside a Wigner chaos of odd order, since random
variables inside such chaoses have all odd moments equal to zero, while
one has, for example, that $\varphi[Z(\lambda)^3] = \lambda$; see
Remark~\ref{calculmoment}(ii).

As a consequence of Theorem~\ref{mainthm}, we will be able to prove the
following result, stating that Wigner chaoses of order greater than 2
do not contain any nonzero Poisson random variable.

\begin{prop}\label{Pnochi2}
Let $q\geq4$ be even, and let $F$ be a nonzero random variable in the
$q$th Wigner chaos. Then $F$ cannot have a free Poisson distribution.
\end{prop}

As pointed out in Remark~\ref{rchi2} below, centered Poisson random
variables with integer rate can be realized as elements of the second
Wigner chaos. As a consequence, Proposition~\ref{Pnochi2} implies that
the second Wigner chaos contains random variables whose distribution
is not shared by any element of higher chaoses. This result parallels
the findings in~\cite{KNPS}, where it is proved that Wigner chaoses of
order $\geq2$ do not contain any nonzero semicircular random
variable. Note that, at the present time, it is not known in general
whether two
nonzero random variables belonging to two distinct Wigner chaoses have
necessarily different laws.

\begin{rem}
We are still far from understanding the exact structure of the free
Wigner chaos. For instance, almost nothing is
known about the regularity of the distributions associated with the
elements of a fixed Wigner chaos.
In particular, we ignore whether such laws may have atoms or are indeed
absolutely continuous (as are those in the classical Wiener chaos).
Further references related to the subject of the present paper are
\cite
{DN,DNNcmp}.
\end{rem}

\subsection{The free probability setting} Our main reference for free
probability is the monograph by Nica and Speicher~\cite{nicaspeicher},
to which the reader is referred for any unexplained notion or result.
We shall also use a notation which is consistent with the one adopted
in~\cite{KNPS}.\vadjust{\goodbreak}

For the rest of the paper, we consider as given a so-called (tracial)
\textit{$W^*$-probability space} $(\mathscr{A},\varphi)$, where:
$\mathscr
{A}$ is a von Neumann algebra of operators (with involution $X\mapsto
X^*$), and $\varphi$ is a unital linear functional on $\mathscr{A}$
with the properties of being \textit{weakly continuous}, \textit{positive}
[i.e., $\varphi(XX^*)\geq0$ for every $X\in\mathscr{A}$], \textit{faithful}
[i.e., such that the relation $\varphi(XX^*) = 0$ implies
$X=0$] and \textit{tracial} [i.e., $\varphi(XY) = \varphi(YX)$, for every
$X,Y\in\mathscr{A}$].

As usual in free probability, we refer to the self-adjoint elements of
$\mathscr{A}$ as \textit{random variables}. Given a random variable $X$ we
write $\mu_X$ to indicate the \textit{law} (or \textit{distribution}) of $X$,
which is defined as the unique Borel probability measure on $\R$ such
that, for every integer $m\geq0$, $\varphi(X^m) = \int_\R x^m\mu_X(dx)$; see,
for example,~\cite{nicaspeicher}, Proposition 3.13.

We say that the unital subalgebras $\mathscr{A}_1,\ldots,\mathscr{A}_n$ of
$\mathscr{A}$ are \textit{freely independent} whenever the following
property holds: let $X_1,\ldots,X_m$ be a finite collection of elements
chosen among the $\mathscr{A}_i$'s in such a way that (for
$j=1,\ldots,m-1$) $X_j$ and $X_{j+1}$ do not come from the same $\mathscr
{A}_i$ and $\varphi(X_j) = 0$ for $j=1,\ldots,m$; then $\varphi(X_1\cdot
\cdot\cdot X_m)=0$. Random variables are said to be freely independent
if they generate freely independent unital subalgebras of $\mathscr{A}$.

\subsection{Plan}

The rest of the paper is organized as follows. In Section \ref
{SPoisson} we provide a characterization of centered free Poisson distributions
in terms of noncrossing partitions. Section~\ref{swigner} deals with
free Brownian motion and Wigner chaos. Section~\ref{sproofs} contains
the proofs of the main results of the paper (i.e., Theorem \ref
{mainthm} and Proposition~\ref{Pnochi2}), whereas Section \ref
{slemmas} is
devoted to some auxiliary lemmas.

\section{Semicircular and centered free Poisson distributions}\label
{SPoisson}

The following definition contains most of the combinatorial objects
that are used throughout the text.

\begin{defi}\label{DC}
(i)
Given an integer $m \geq1$, we write $[m] = \{1,\ldots,m\}$. A~\textit{partition}
of $[m]$ is a collection of nonempty and disjoint subsets
of $[m]$, called \textit{blocks}, such that their union is equal to $[m]$.
The cardinality of a block is called \textit{size}. A block is said to be
a \textit{singleton} if it has size one.\vspace*{-6pt}
\begin{longlist}[(iii)]
\item[(ii)] A partition $\pi$ of $[m]$ is said to be \textit{noncrossing}
if one cannot find integers $p_1,q_1,p_2,q_2$ such that: (a) $1\leq p_1
< q_1 < p_2 < q_2\leq m$, (b) $p_1,p_2$ are in the same block of $\pi$,
(c) $q_1,q_2$ are in the same block of $\pi$ and (d) the $p_i$'s are
not in the same block of $\pi$ as the $q_i$'s. The collection of the
noncrossing partitions of $[m]$ is denoted by $NC(m)$, $m\geq1$.

\item[(iii)] For every $m\geq1$, the quantity $C_m = | NC(m) |$, where
$ |A |$ indicates the cardinality of a set $A$, is called the $m$th
\textit{Catalan number}. One sets by convention $C_0 = 1$. Also, recall
the explicit expression $C_m = \frac{1}{m+1} \bigl({2m\atop m} \bigr)$.\vspace*{1pt}

\item[(iv)] We define the sequence $\{R_m \dvtx m\geq0\}$ as follows: $R_0
= 1$, and, for \mbox{$m\geq1$}, $R_m$ is equal to the number of partitions in
$NC(m)$ having no singletons.\vadjust{\goodbreak}

\item[(v)] For every $m\geq1$ and every $j=1,\ldots,m$, we define
$R_{m,j}$ to be the number of noncrossing partitions of $[m]$ with
exactly $j$ blocks and with no singletons. Plainly, $R_m = \sum_{j=1}^m
R_{m,j}$. Also, when $m$ is even, one has that $R_{m,j} = 0$ for every
$j>m/2$; when $m$ is odd, then $R_{m,j} = 0$ for every
$j>(m-1)/2$.\vspace*{-2pt}
\end{longlist}
\end{defi}

\begin{example}\label{suf}
One has that:
\begin{longlist}[--]
\item[--] $R_1 = R_{1,1} = 0$, since $\{\{1\}\}$ is the only partition
of $[1]$, and such a partition is composed of exactly one singleton;
\item[--] $R_2 = R_{2,1}= 1$, since the only partition of $[2]$ with no
singletons is $\{\{1,2\}\}$;
\item[--] $R_3 =R_{3,1} = 1$, since the only partition of $[3]$ with no
singletons is $\{\{1,2,3\}\}$;
\item[--] $R_4 = 3$, since the only noncrossing partitions of $[4]$
with no singletons are $\{\{1,2,3,4\}\}$, $\{\{1,2\},\{3,4\}\}$ and $\{
\{1,4\},\{2,3\}\}$. This implies that $R_{4,1} = 1$ and $R_{4,2} =
2$.\vspace*{-2pt}
\end{longlist}
\end{example}

The integers $\{R_m \dvtx m\geq0\}$ are customarily called the
\textit{Riordan numbers}. A detailed analysis of the combinatorial
properties of Riordan numbers is provided in the paper by Bernhart
\cite {bernhart}; however, it is worth noting that the discussion to
follow is self-contained, in the sense that no previous knowledge of
the combinatorial properties of the sequence $\{R_m\}$ is required.

Given a random variable $X$, we denote by $\{\kappa_m(X) \dvtx m\geq1\}$
the sequence of the \textit{free cumulants} of $X$. We recall
(see~\cite{nicaspeicher}, page 175) that the free cumulants of $X$ are
completely determined by the following relation: for every $m\geq1$,
%
\begin{equation}
\label{emomcum} \varphi\bigl(X^m\bigr) = \sum
_{\pi= \{b_1,\ldots,b_j\}\in NC(m)} \prod_{i=1}^j
\kappa_{|b_i|}(X),
\end{equation}
where $|b_i|$ indicates the size of the block $b_i$ of the noncrossing
partition $\pi$. It is clear from (\ref{emomcum}) that the sequence
$\{
\kappa_m(X) \dvtx m\geq1\}$ completely determines the moments of $X$ (and
vice-versa).\vspace*{-2pt}

\begin{defi}\label{DZ}
(i) The centered \textit{semicircular distribution} of parameter
\mbox{$t>0$}, denoted by $S(0,t)(dx)$, is the probability distribution given
by
\[
S(0,t) (dx) = (2\pi t)^{-1}\sqrt{4t-x^2}\,dx,\qquad |x|< 2
\sqrt{t}.
\]
We recall the following classical relation:
\[
\int_{-2\sqrt{t}}^{2\sqrt{t}} x^{2m} S(0,t) (dx) =
C_m t^m,
\]
where $C_m$ is the $m$th Catalan number [so that, e.g., the second
moment of $S(0,t)$ is $t$]. Since the odd moments of $S(0,t)$ are all
zero, one deduces from the previous relation and (\ref{emomcum})
(e.g., by recursion) that the free cumulants of a random variable~$s$
with law $S(0,t)$ are all zero, except for $\kappa_2(s) =\varphi
(s^2)= t$.\vadjust{\goodbreak}

\begin{longlist}
\item[(ii)] The \textit{free Poisson distribution} with rate $\lambda>0
$, denoted by $P(\lambda)(dx)$ is the probability distribution defined
as follows: (a) if $\lambda\in(0,1]$, then $P(\lambda) = (1-\lambda
)\delta_0 + \lambda\widetilde{\nu}$, and (b) if $\lambda> 1$, then
$P(\lambda) = \widetilde{\nu}$, where $\delta_0$ stands for the Dirac
mass at $0$. Here, $\widetilde{\nu}(dx) =
(2\pi x)^{-1}\sqrt{4\lambda-(x-1-\lambda)^2}\,dx,   x\in
((1-\sqrt {\lambda})^2,(1+\sqrt{\lambda})^2 )$. If $X_\lambda$ has the
$P(\lambda)$ distribution, then~\cite{nicaspeicher}, Proposition 12.11,
implies that
%
\begin{equation}
\label{ecumPOISS} \kappa_m(X_\lambda) = \lambda, \qquad m\geq1.
\end{equation}
From now on, we will denote by $Z(\lambda)$ a random variable having
the law of $X_\lambda- \lambda1$
(centered free Poisson distribution), where $1$ is the unity of
$\mathscr{A}$. Plainly, $\kappa_1[Z(\lambda)] = \varphi[Z(\lambda
)] =
0$, and $\kappa_2[Z(\lambda)] = \varphi[X_\lambda^2] - \lambda^2$ is
the variance of $X_\lambda$.
\end{longlist}
\end{defi}

Note that both $S(0,t)$ and $P(\lambda)$ are compactly supported, and
therefore are uniquely determined by their moments
(by the Weierstrass theorem). Definition~\ref{DZ}(ii) is taken from
\cite{nicaspeicher}, Definition 12.12. As discussed in the
\hyperref[sec1]{Introduction}, the choice of the denomination
``free Poisson'' comes from the following two facts: (1) $P(\lambda)$
can be obtained as the limit of the free convolution of Bernoulli
distributions (see~\cite{nicaspeicher}, Proposition 12.11), and (2) the
classical Poisson distribution of parameter $\lambda$ has (classical)
cumulants all equal to $\lambda$ (see, e.g.,~\cite{PecTaq}, Section 3.3).
As already pointed out, the free Poisson law is also called the
``Marchenko--Pastur distribution.''

The following statement contains a characterization of the moments of
$Z(\lambda)$, and shows that, when $\lambda$ is integer, then
$Z(\lambda)$ is the free equivalent of a classical centered $\chi^2$
random variable with $\lambda$ degrees of freedom. This last fact could
alternatively be deduced from~\cite{nicaspeicher}, Proposition 12.13,
but here we prefer to provide a self-contained argument.

\begin{prop}\label{Po} Let the notation of Definitions~\ref{DC} and~\ref{DZ} prevail. Then, for every real $\lambda>0$ and every
integer $m\geq1$,
%
\begin{equation}
\label{eMP} \varphi\bigl[Z(\lambda)^m\bigr] = \sum
_{j=1}^m \lambda^j R_{m,j}.
\end{equation}
Let $p=1,2,\ldots$ be an integer; then $Z(p)$ has the same law as $\sum_{i=1}^p(s_i^2 -1)$, where $s_1,\ldots,s_p$ are $p$ freely independent
random variables with the $S(0,1)$ distribution, and $1$ is the unit of
$\mathscr{A}$.
\end{prop}
\begin{pf} From (\ref{ecumPOISS}), one deduces that $\kappa_m[Z(\lambda)] = \lambda$ for every $m\geq2$. Since $\kappa_1[Z(\lambda
)] =0$,
we infer from (\ref{emomcum}) that
\[
\varphi\bigl[Z(\lambda)^m\bigr] =\sum_{j=1}^m
\sum_{\pi= \{b_1,\ldots,b_j\}
\in
NC(m)} \lambda^j \mathbf{1}_{\{\pi\ \mathrm{has\ no\ singletons\}}},
\]
which immediately yields (\ref{eMP}). To prove the last part of the statement,
consider first the case $p=1$, write $s = s_1$ and fix an integer
$m\geq2$. In order to build a noncrossing partition of $[m]$, say
$\pi
$, one has to perform the following three steps: (a) choose an integer
$j\in\{0,\ldots,m\}$, denoting the number of singletons of $\pi$, (b)~choose the $j$ singletons of $\pi$ among the $m$ available integers
[this can be done in exactly $\bigl({m\atop j}\bigr)$ distinct ways], (c) build a
noncrossing partition of the remaining $m-j$ integers with blocks at
least of size~$2$ (this can be done in exactly $R_{m-j}$ distinct
ways). Since $C_0 = R_0 =1$ and $C_1 =1 = R_0+R_1$, it follows that
Catalan and Riordan numbers are linked by the following relation: for
every $m\geq0$,
%
\begin{equation}
\label{u} C_m = \sum_{j=0}^m
\pmatrix{m
\cr
j} R_{m-j} = \sum_{j=0}^m
\pmatrix {m
\cr
j} R_{j},
\end{equation}
where the last equality follows from $\bigl({m\atop j}\bigr) = \bigl({m\atop m-j}\bigr)$.
By inversion, one therefore deduces that
\[
R_m = \sum_{j=0}^m \pmatrix{m
\cr
j} (-1)^{m-j}C_{j},\qquad m\geq0.
\]
Therefore
\begin{eqnarray*}
\varphi\bigl[\bigl(s^2-1\bigr)^m\bigr]&=& \sum
_{j=0}^m \pmatrix{m
\cr
j}(-1)^{m-j}\varphi
\bigl(s^{2j}\bigr)=\sum_{j=0}^m
\pmatrix{m
\cr
j} (-1)^{m-j}C_j
\\
&=&R_m=\sum_{j=1}^m
R_{m,j} = \varphi\bigl[Z(1)^m\bigr],
\end{eqnarray*}
from which we infer that $s^2-1\stackrel{\mathrm{law}}{=}Z(1)$, yielding the
desired conclusion when $p=1$.
Let us now consider the general case, that is, $p\geq2$. First recall
that the $m$th free cumulant of the sum of $p$ freely independent
random variables is the sum of the corresponding $m$th cumulants (this
is a consequence of the multilinearity of free cumulants, as well as of
the characterization of free independence in terms of vanishing mixed
cumulants; see~\cite{nicaspeicher}, Theorem~11.16). It follows that, for
any $m\geq2$,
\[
\kappa_m \Biggl(\sum_{i=1}^p
\bigl(s_i^2-1\bigr) \Biggr) = p\times
\kappa_m\bigl(s_1^2-1\bigr)=p\times
\kappa_m\bigl(Z(1)\bigr) = p=\kappa_m\bigl(Z(p)\bigr).
\]
This implies that $\sum_{i=1}^p s_i^2-1\stackrel{\mathrm{law}}{=}Z(p)$, and
the proof of Proposition~\ref{Po} is concluded.
\end{pf}

\begin{rem}\label{calculmoment}
(i) Relation (\ref{u}) is well known; see, for example,
\cite{bernhart},
Section~5, for an alternate proof based on ``difference
triangles.'' Our proof of the relation $R_m = \varphi[Z(1)^m]$ seems to
be new.
\begin{longlist}
\item[(ii)] Using the last two points of Example~\ref{suf}, we deduce
from (\ref{eMP}) that $\varphi[Z(\lambda)^3] =
\lambda R_{3,1}=\lambda$, while $\varphi[Z(\lambda)^4] =\lambda
R_{4,1} + \lambda^2 R_{4,2}=\lambda+2\lambda^2$.
\end{longlist}
\end{rem}

\section{Free Brownian motion and Wigner chaos}\label{swigner}
Our main reference for the content of this section is the paper by
Biane and Speicher~\cite{bianespeicher}.

\begin{defi}[($L^p$ spaces)]\label{dspaces}
(i) For $1\leq p \leq\infty$, we write $L^p(\mathscr
{A},\varphi)$ to indicate the $L^p$ space obtained as the completion of
$\mathscr{A}$ with respect to the norm $\| a\|_p = \varphi
(|a|^p)^{1/p}$, where $ |a|= \sqrt{a^\ast a}$, and $\|\cdot\|_\infty$
stands for the operator norm.\vspace*{-6pt}
\begin{longlist}[(iii)]
\item[(ii)] For every integer $q\geq2$, the space $L^2(\mathbb
{R}_+^q)$ is the collection of all complex-valued functions on $\mathbb
{R}_+^q$ that are square-integrable with respect to the Lebesgue
measure. Given $f\in L^2(\mathbb{R}_+^q)$, we write
\[
f^*(t_1,t_2,\ldots,t_q) =
\overline{f(t_q,\ldots,t_2,t_1)},
\]
and we call $f^*$ the \textit{adjoint} of $f$. We say that an element of
$L^2(\mathbb{R}_+^q)$ is \textit{mirror symmetric} if
\[
f(t_1,\ldots,t_q) = f^*(t_1,
\ldots,t_q)
\]
for almost every vector $(t_1,\ldots,t_q)\in\mathbb{R}_+^q$. Notice that
mirror symmetric functions constitute a Hilbert subspace of
$L^2(\mathbb
{R}_+^q)$.
\item[(iii)] Given $f\in L^2(\mathbb{R}_+^q)$ and $g\in L^2(\mathbb
{R}_+^p)$, for every $r = 1,\ldots,\min(q,p)$ we define the $r$th
\textit{contraction} of $f$ and $g$ as the element of $L^2(\mathbb
{R}_+^{p+q-2r})$ given by
%
\begin{eqnarray}
\label{econtraction2} &&f\cont{r} g(t_1,\ldots,t_{p+q-2r})
\nonumber\\
&& \qquad= \int_{\mathbb{R}_+^r} f(t_1,\ldots,t_{q-r},y_r,y_{r-1},
\ldots,y_1)\\
&&\hspace*{49pt}{}\times g(y_1,y_2,\ldots,y_r,
t_{q-r+1},t_{p+q-2r})\,dy_1\cdots\,dy_r.\nonumber
\end{eqnarray}
One also writes $f\cont{0} g (t_1,\ldots,t_{p+q}) = f\otimes g
(t_1,\ldots,t_{p+q}) = f(t_1,\ldots, t_q)g(t_{q+1},\break  \ldots,t_{p+q})$. In the
following, we shall use the notation $f\cont{0} g$ and $f\otimes g$
interchangeably. Observe that, if $p=q$, then $f\cont{p} g = \langle
f,g^{*}\rangle_{L^2(\R_+^q)}$.
\end{longlist}
\end{defi}

A \textit{free Brownian motion} $S$ on $(\mathscr{A},\varphi)$ consists
of: (i) a filtration $\{\mathscr{A}_t \dvtx t\geq0\}$ of von Neumann
sub-algebras of $\mathscr{A}$ (in particular, $\mathscr{A}_u \subset
\mathscr{A}_t$, for $0\leq u<t$), (ii)~a~collection $S = \{S(t) \dvtx t\geq
0\}$ of self-adjoint operators such that:
\begin{longlist}[--]
\item[--] $S(t)\in\mathscr{A}_t$ for every $t$;
\item[--] for every $t$, $S(t)$ has a semicircular distribution
$S(0,t)$, with mean zero and variance $t$;
\item[--] for every $0\leq u<t$, the ``increment'' $S(t) - S(u)$ is
freely independent of $\mathscr{A}_u$, and has a semicircular
distribution $S(0,t-u)$, with mean zero and variance $t-u$.
\end{longlist}

For every integer $q\geq1$, the collection of all random variables of
the type $I^S_q(f) = I(f)$, $f \in L^2(\mathbb{R}_+^q)$, is called the
$q$th \textit{Wigner chaos} associated with~$S$, and is defined according
to~\cite{bianespeicher}, Section 5.3, namely:
\begin{longlist}[--]
\item[--] first define $I(f) = (S({b_1}) - S({a_1}))\cdots(S({b_q}) -
S({a_q}))$, for every function $f$ having the form
%
\begin{equation}
\label{esimple} f(t_1,\ldots,t_q) = \prod
_{i=1}^q \mathbf{1}_{(a_i,b_i)}(t_i),
\end{equation}
where the intervals $(a_i,b_i)$, $i=1,\ldots,q$, are pairwise disjoint;
\item[--] extend linearly the definition of $I(f)$ to ``simple
functions vanishing on diagonals,'' that is, to functions $f$ that are
finite linear combinations of indicators of the type (\ref{esimple});
\item[--] exploit the isometric relation
%
\begin{eqnarray}
\label{efreeisometry} \bigl\langle I(f_1),I(f_2) \bigr
\rangle_{L^2(\mathscr{A},\varphi)}&=& \varphi \bigl(I(f_1)^*I(f_2)
\bigr)= \varphi \bigl(I\bigl(f_1^*\bigr)I(f_2) \bigr)
\nonumber
\\[-8pt]
\\[-8pt]
\nonumber
&=&
\langle f_1,f_2 \rangle_{L^2(\mathbb{R}_+^q)},
\end{eqnarray}
where $f_1,f_2$ are simple functions vanishing on diagonals, and use a
density argument to define $I(f)$ for a general $f\in L^2(\mathbb{R}_+^q)$.
\end{longlist}

Observe that relation (\ref{efreeisometry}) continues to hold for
every pair $f_1,f_2 \in L^2(\mathbb{R}_+^q)$. Moreover, the above
sketched construction implies that $I(f)$ is self-adjoint if and only
if $f$ is mirror symmetric. Finally, we recall the following
fundamental multiplication formula, proved in~\cite{bianespeicher}. For
every $f\in L^2(\R_+^q)$ and $g\in L^2(\R_+^p)$, where $q,p\geq1$,
%
\begin{equation}
\label{emult} I(f)I(g) = \sum_{r=0}^{\min(q,p)} I(f
\cont{r}g).
\end{equation}

\begin{rem}\label{rchi2}Let $\{e_i \dvtx 1,\ldots,p\}$ be an orthonormal
system in $L^2(\R_+)$. Then, the random variables $s_i=I(e_i)$,
$i=1,\ldots,p$, have the $S(0,1)$ distribution and are freely independent.
Moreover, the product formula (\ref{emult}) implies that
\[
\sum_{i=1}^p\bigl(s_i^2-1
\bigr) = I \Biggl(\sum_{i=1}^p
e_i\otimes e_i \Biggr),
\]
and therefore that the double integral $I (\sum_{i=1}^p
e_i\otimes
e_i )$ has a centered free Poisson distribution with rate $p$.
\end{rem}

\section{Proof of the main results}\label{sproofs}
\subsection{\texorpdfstring{Proof of Theorem \protect\ref{mainthm}}{Proof of Theorem 1.4}}

In the free probability setting (see, e.g.,~\cite{nicaspeicher}, Definition
8.1) convergence in distribution is equivalent to the
convergence of moments, so that $I(f_n)$ converges in distribution to
$Z(\lambda)$ if and only if $\varphi(I(f_n)^m) \to\varphi(Z(\lambda
)^m)$, for every $m\geq1$. In particular, convergence in distribution
implies \mbox{$\varphi(I(f_n)^4)- 2\varphi(I(f_n)^3) \to\varphi(Z(\lambda
)^4)- 2\varphi(Z(\lambda)^3) = 2\lambda^2-\lambda$}.

Now assume that $\varphi[I(f_n)^4]-2\varphi[I(f_n)^3]\to2\lambda^2-\lambda$.
We have to show that, for every $m\geq2$, $\varphi[I(f_n)^m] \to
\varphi[Z(\lambda)^m]$. Iterative applications of the product formula
(\ref{emult}) yield
%
\begin{equation}
\label{ae} I(f_n)^m=\sum_{(r_1,\ldots,r_{m-1})\in A_m}
I \bigl(\bigl(\cdots\bigl((f_n\cont {r_1}f_n)
\cont{r_2}f_n\bigr)\cdots f_n\bigr)
\cont{r_{m-1}}f_n \bigr),
\end{equation}
where
\begin{eqnarray*}
A_m&=& \bigl\{ (r_1,\ldots,r_{m-1})\in\{0,1,
\ldots,q\}^{m-1}\dvtx r_2\leq2q-2r_1,\\
&&\hspace*{7pt}r_3\leq3q-2r_1-2r_2,\ldots, r_{m-1}\leq(m-1)q-2r_1-\cdots-2r_{m-2}
\bigr\};
\end{eqnarray*}
note that (\ref{ae}) was proved in~\cite{KNPS}, formula (1.10). We
therefore deduce that
\[
\varphi\bigl[I(f_n)^m\bigr]=\sum
_{(r_1,\ldots,r_{m-1})\in B_m} \bigl(\cdots\bigl((f_n\cont{r_1}f_n)
\cont{r_2}f_n\bigr)\cdots f_n\bigr)
\cont{r_{m-1}}f_n,
\]
with $B_m= \{(r_1,\ldots,r_{m-1})\in A_m\dvtx 2r_1+\cdots
+2r_{m-1}=mq \}$. The previous equality is a consequence of the
following fact: in the sum on the right-hand side of~(\ref{ae}), the
elements indexed by $B_m$ correspond to constants, whereas the elements
indexed by $A_m\setminus B_m$ are genuine multiple Wigner integrals,
and therefore have $\varphi$-expectation equal to zero.
We further decompose $B_m$ as follows: $B_m=D_m\cup E_m$, with
$D_m=B_m\cap\{0,\frac{q}2,q\}^{m-1}$ and $E_m=B_m\setminus D_m$,
so that
%
\begin{eqnarray}
\label{doublesum}
\varphi\bigl[I(f_n)^m
\bigr]&=& \sum_{(r_1,\ldots,r_{m-1})\in D_m} \bigl(\cdots\bigl((f_n
\cont{r_1}f_n)\cont{r_2}f_n
\bigr)\cdots f_n\bigr)\cont {r_{m-1}}f_n
\nonumber
\\[-8pt]
\\[-8pt]
\nonumber
&&{}+\sum_{(r_1,\ldots,r_{m-1})\in E_m} \bigl(\cdots\bigl((f_n
\cont{r_1}f_n)\cont{r_2}f_n
\bigr)\cdots f_n\bigr)\cont {r_{m-1}}f_n.
\end{eqnarray}
By the forthcoming Lemma~\ref{lm2}, we have
$\|f_n\cont{q/2} f_n -f_n\|\to0$
and $\|f_n\cont{r} f_n\|\to0$ for $r\in\{1,\ldots,q-1\}\setminus\{
\frac
{q}2\}$.
The conclusion is then obtained by observing that the first sum in
(\ref
{doublesum}) converges to $\varphi[Z(\lambda)^m]$ by
Proposition~\ref{Po} and the forthcoming Lemma~\ref{lm4},
whereas
the second sum converges to zero by the forthcoming Lemma~\ref{lm5}.

\subsection{\texorpdfstring{Proof of Proposition \protect\ref{Pnochi2}}{Proof of Proposition 1.5}}

Assume that $F = I(f)$, where $f$ is a mirror symmetric element of
$L^2(\R_+^q)$ for some even $q\geq4$, and also that $\varphi[F^2] =
\|f\|^2_{L^2(\R_+^q)} = \lambda>0$. If $F$ had the same\vspace*{-3pt} law of
$Z(\lambda)$, then $\varphi(F^4) - 2\varphi(F^3) = 2\lambda^2-\lambda$, and
the forthcoming Lemma~\ref{lm2} would imply that $\|f\cont{q/2}f-f\|_{L^2(\R^q_+)} = 0$ and $\|f \cont{r} f\|_{L^2(\R^{2q-2r}_+)} = 0$ for
all $r\in\{1,\ldots,q-1\}\setminus\{ \frac{q}2\}$. As shown\vspace*{-2pt} in
\cite{KNPS},
Proof of Corollary 1.7, the relation $\|f \cont{q-1} f\|_{L^2(\R
^{2}_+)} = 0$ implies that necessarily $f = 0$, and therefore that
$F=0$. Since $\varphi[F^2] = \lambda>0$ we have achieved a
contradiction, and the proof is complete.

\section{Ancillary lemmas}\label{slemmas}

This section collects some technical results that are used in the proof
of Theorem~\ref{mainthm}.

\begin{lemma}\label{lm2}
Let $q\geq2$ be an even integer, and consider a sequence $\{f_n \dvtx n\geq1\} \subset L^2(\R^q_+)$ of mirror symmetric
functions such that $\|f_n\|^2_{L^2(\R_+^q) }= \lambda>0$ for every
$n$. As $n\to\infty$, one has that
\[
\varphi\bigl[I(f_n)^4\bigr]-2 \varphi
\bigl[I(f_n)^3\bigr] \to2\lambda^2-\lambda
\]
if and only if $\|f_n\cont{q/2}f_n-f_n\|_{L^2(\R^q_+)}\to0$ and $\|f_n
\cont{r} f_n\|_{L^2(\R^{2q-2r}_+)}\to0$ for
all $r\in\{1,\ldots,q-1\}\setminus\{ \frac{q}2\}$.
\end{lemma}
\begin{pf}The product formula yields
\[
I(f_n)^2 - I(f_n)= \lambda+
I(f_n\cont{0} f_n)+ I(f_n\cont
{q/2}f_n-f_n) + \mathop{\sum
_{1\leq r\leq q-1}}_{r\neq q/2} I(f_n\cont{r}
f_n).
\]
Using the isometry property and the fact that multiple Wigner integrals
of different orders are orthogonal in $L^2(\mathscr{A},\varphi)$, we
deduce that
\begin{eqnarray*}
\label{fourthpower} &&\varphi\bigl[\bigl(I(f_n)^2-I(f_n)
\bigr)^2\bigr]
\\
&&\qquad= \lambda^2 + \|f_n\cont{0} f_n
\|_{L^2(\R^{2q}_+)}^2 + \|f_n\cont {q/2}f_n-f_n
\|_{L^2(\R^q_+)}^2\\
&&\qquad\quad{}+ \mathop{\sum_{1\leq r\leq q-1}}_{r\neq q/2}
\|f_n\cont{r} f_n\|_{L^2(\R^{2q-2r}_+)}^2
\\
&&\qquad= 2\lambda^2+\|f_n\cont{q/2}f_n-f_n
\|_{L^2(\R^{q}_+)}^2+ \mathop{\sum_{1\leq r\leq q-1}}_{r\neq q/2}
\|f_n \cont{r} f_n\|_{L^2(\R^{2q-2r}_+)}^2,
\end{eqnarray*}
and the desired conclusion follows because $\varphi[I(f_n)^2]=\|f_n\|^2_{L^2(\R_+^q) }=\lambda$.
\end{pf}

\begin{lemma}\label{lm4}
Let $m\geq2$ be an integer, let $q\geq2$ be an even integer and
recall the notation adopted in (\ref{doublesum}). Assume
$\{f_n \dvtx n\geq1\} \subset L^2(\R_+^q)$ is a sequence of mirror
symmetric functions satisfying $\|f_n\|^2_{L^2(\R_+^q) }= \lambda>0$
for every $n$.\vspace*{-2pt} If
$\|f_n\cont{q/2}f_n-f_n\|^2_{L^2(\R_+^q) }\to0$ as $n\to\infty$,
then
%
\begin{eqnarray}
\label{jessold}&& \sum_{(r_1,\ldots,r_{m-1})\in D_m} \bigl(\cdots
\bigl((f_n\cont {r_1}f_n)\cont
{r_2}f_n\bigr)\cdots f_n\bigr)
\cont{r_{m-1}}f_n\to\varphi\bigl[Z(\lambda)^m
\bigr]
\nonumber
\\[-9pt]
\\[-9pt]
\nonumber
&&\qquad= \sum_{j=1}^m \lambda^j
R_{m,j}
\end{eqnarray}
as $n\to\infty$.
\end{lemma}
\begin{pf} Assume that $f_n\cont{q/2}f_n \approx f_n$ (given two
sequences $\{a_n\}$ and $\{b_n\}$ with values in some normed vector
space, we write
$a_n\approx b_n$ to indicate that $a_n-b_n\to0$ with respect to the
associated norm), and consider
$(r_1,\ldots,\break r_{m-1})\in D_m$. We now claim that
%
\begin{equation}
\label{eas} \bigl(\cdots\bigl((f_n\cont{r_1}f_n)
\cont{r_2}f_n\bigr)\cdots f_n\bigr)
\cont{r_{m-1}}f_n \to\lambda^j,
\end{equation}
where $j$ equals the number of the entries of $(r_1,\ldots,r_{m-1}) $
that are equal to $q$.
To see why (\ref{eas}) holds, write $G_0 := f_n$, $G_1 := f_n\cont
{r_1}f_n, \ldots, G_{m-1}:=(\cdots((f_n\cont{r_1}f_n)\cont
{r_2}f_n)\cdots f_n)\cont{r_{m-1}}f_n$, and observe that the following
facts take place:
\begin{longlist}[(iii)]
\item[(i)] Since $r_j \in\{0, q/2, q\}$, every function $G_j$ is
either a constant or a multiple of an object of the type $H_1 \otimes
\cdots\otimes H_l$, where $l\geq1$,  and every $H_i$ ($i=1,\ldots,l$) is
either equal to $f_n$ or to an iterated contraction of the type
\[
\underbrace{f_n\cont{q/2}\cdots\cont{q/2} f_n
}_{k\ \mathrm{contractions}}
\]
for some $k\geq1$. In particular, every $H_i$ is a function in $q$ variables.

\item[(ii)] If $G_j = c$ is a constant, then necessarily $r_{j+1} = 0$
and $G_{j+1} = c\times f_n$.

\item[(iii)] If $c$ is a constant, $G_j = c\times H_1 \otimes\cdots
\otimes H_l$ and $r_{j+1} = q$, then $G_{j+1} = c\langle H_l, f_n
\rangle_{L^2(\R^q_+)}\times H_1 \otimes\cdots\otimes H_{l-1}$.

\item[(iv)] If $c$ is a constant, $G_j = c\times H_1 \otimes\cdots
\otimes H_l$ and $r_{j+1} = q/2$, then $G_{j+1} = c H_1 \otimes\cdots
\otimes(H_{l}\cont{q/2} f_n ) $.

\item[(v)] since $(r_1,\ldots,r_m)\in B_m$, the quantity $G_{m-1}$ is
necessarily a constant given by the product of $j$ scalar products
having either the form $\langle f_n, f_n \rangle_{L^2(\R^q_+)} =
\lambda
$, or
%
\begin{equation}
\label{bs} \bigl\langle \underbrace{f_n\cont{q/2}\cdots\cont{q/2}
f_n }_{k\ \mathrm{contractions}} , f_n \bigr\rangle_{L^2(\R^q_+)}
\end{equation}
for some $k\geq1$.\vadjust{\goodbreak}
\end{longlist}
Using the two relations $f_n\cont{q} f_n = \|f_n\|_{L^2(\R
_+^q)}^2=\lambda$ and $f_n\cont{q/2}f_n \approx f_n$, one sees
immediately that the left-hand side of (\ref{bs}) converges to
$\lambda
$, as $n\to\infty$, so that relation (\ref{eas}) is proved.

As a consequence, for every $m\geq2$, there exists a polynomial
$w_m(\lambda)$ (independent of $q$) such that, for every sequence $\{
f_n\}$ as in the statement,
\[
\sum_{(r_1,\ldots,r_{m-1})\in D_m} \bigl(\cdots
\bigl((f_n\cont {r_1}f_n)\cont
{r_2}f_n\bigr)\cdots f_n\bigr)
\cont{r_{m-1}}f_n\to w_m(\lambda).
\]
Now consider the case $q=2$ and $f_n = f = \sum_{i=1}^p e_i\otimes
e_i$, where $p\geq1$ and $\{e_i \dvtx i=1,\ldots,p\}$ is an orthonormal
system in
$ L^2(\R_+^q)$. The following three facts take place: (a)
$I(\sum_{i=1}^p e_i\otimes e_i)$ has the same law as $Z(p)$ (see Remark~\ref{rchi2}), (b) $\|f\|^2_{L^2(\R_+^2)} = p$ and (c) $f\cont{1}f =f$.
Since $E_m=\varnothing$ for $q=2$, the previous discussion [combined with
(\ref{doublesum}) and Proposition~\ref{Po}] yields that, for every
$m\geq2$, $w_m(p) = \varphi[Z(p)^m] = \sum_{j=1}^{m} p^j R_{m,j}$, for
every $p=1,2,\ldots.$ Since two polynomials coinciding on a countable set
are necessarily equal, we deduce that $w_m(\lambda) = \varphi
[Z(\lambda
)^m]$ for every $\lambda>0$.
\end{pf}

\begin{rem}By inspection of the arguments used in the proof of
Lemma~\ref{lm4}, one deduces that $R_m = |D_m|$, for every $m\geq2$.
\end{rem}

\begin{lemma}
\label{lm5}
Let $m\geq2$ be an integer, let $q\geq2$ be an even integer and
recall the notation adopted in (\ref{doublesum}). Assume $\{f_n \dvtx n\geq
1\} \subset L^2(\R_+^q)$ is a sequence of mirror symmetric functions
satisfying $\|f_n\|^2_{L^2(\R_+^q) }= \lambda>0$ for every $n$. If
$(r_1,\ldots,r_{m-1})\in E_m$ and if $\|f_n\cont{r} f_n\|_{L^2(\R
_+^{2q-2r}) }\to0$ for all $r\in\{1,\ldots,q-1\}\setminus\{\frac
{q}2\}
$, then
\[
\bigl(\cdots\bigl((f_n\cont{r_1}f_n)
\cont{r_2}f_n\bigr)\cdots f_n\bigr)
\cont{r_{m-1}}f_n \to0 \qquad\mbox{as $n\to\infty$}.
\]
\end{lemma}
\begin{pf}
This lemma is a straightforward extension of~\cite{KNPS}, Proposition~2.2.
Indeed, going back to the definition of $E_m$ and using the language
introduced in~\cite{KNPS}, observe first that one can rewrite the quantity
\[
\bigl(\cdots\bigl((f_n\cont{r_1}f_n)
\cont{r_2}f_n\bigr)\cdots f_n\bigr)
\cont{r_{m-1}}f_n
\]
as
\[
\int_\pi f_n^{\otimes m},
\]
where $\pi$ is a (uniquely defined) noncrossing pairing such that: (i)
$\pi$ respects $q^{\otimes m}$;
and (ii) $\pi$ is such that there exists two blocks of $q^{\otimes m}$
that are linked
by $r$ pairs for some $r\in\{1,\ldots,q-1\}\setminus\{\frac{q}2\}$.
The desired conclusion then follows by adapting the proof of
\cite{KNPS},
Proposition 2.2,
to this slightly different context.
\end{pf}

\section*{Acknowledgement}
We are grateful to an anonymous referee for a careful reading and a
number of helpful suggestions.

%

%


\printaddresses

\end{document}